\documentclass[final,1pt,times]{elsarticle}

\usepackage{epsfig, graphicx, subfigure}
\usepackage{verbatim, setspace}
\usepackage{amsmath,amssymb}
\usepackage[colorlinks=true, breaklinks=true, backref=page]{hyperref}

\newcommand{\R}{{\mathbb R}}

\newcommand{\T}{{\mathbb T}}

\newcommand{\ic}{{\mathrm i}}
\newcommand{\dd}{{\mathrm d}}

\newtheorem{theorem}{Theorem}

\newtheorem{lemma}[theorem]{Lemma}

\newproof{proof}{Proof}

\journal{}

\begin{document}

\begin{frontmatter}



\title{An asymptotic expansion for the expected number of real zeros of real random polynomials spanned by OPUC}


\author[add]{Hanan Aljubran}

\ead{haljubra@umail.iu.edu}

\author[add]{Maxim L. Yattselev\corref{cor}}

\address[add]{Department of Mathematical Sciences, Indiana University-Purdue University Indianapolis, 402~North Blackford Street, Indianapolis, IN 46202, USA}

\ead{maxyatts@iupui.edu}

\cortext[cor]{Corresponding author.}

\begin{abstract}
Let \( \{\varphi_i\}_{i=0}^\infty \) be a sequence of orthonormal polynomials on the unit circle with respect to a positive Borel measure \( \mu \) that is symmetric with respect to conjugation. We study asymptotic behavior of the expected number of real zeros, say \( \mathbb E_n(\mu) \), of random polynomials
\[
P_n(z) := \sum_{i=0}^n\eta_i\varphi_i(z),
\]
where \( \eta_0,\dots,\eta_n \) are i.i.d. standard Gaussian random variables.  When \( \mu \) is the acrlength measure such polynomials are called Kac polynomials and it was shown by Wilkins that \( \mathbb E_n(|\dd\xi|) \) admits an asymptotic expansion of the form
\[
\mathbb E_n(|\dd\xi|) \sim \frac2\pi\log(n+1) + \sum_{p=0}^\infty A_p(n+1)^{-p}
\]
(Kac himself obtained the leading term of this expansion). In this work we generalize the result of Wilkins to the case where \( \mu \) is absolutely continuous with respect to arclength measure and its Radon-Nikodym derivative extends to a holomorphic non-vanishing function in some neighborhood of the unit circle. In this case \( \mathbb E_n(\mu) \) admits an analogous expansion with coefficients the \( A_p \) depending on the measure \( \mu \) for \( p\geq 1 \) (the leading order term and \( A_0 \) remain the same).
\end{abstract}


\end{frontmatter}


\section{Introduction and Main Results}
\label{}

Random polynomials is a relatively old subject with initial contributions by by Bloch and P\'olya, Littlewood and Offord, Erd\"os and Offord, Arnold, Kac, and many other authors. An interested reader can find a well referenced early history of the subject in the books by Bharucha-Reid and Sambandham \cite{Bharucha-ReidSambandham}, and by Farahmand \cite{Farahmand}. In \cite{Kac43}, Kac considered random polynomials of the form
\[
P_n(z) = \eta_0 + \eta_1 z + \cdots +\eta_nz^n,
\]
where \( \eta_i \) are i.i.d. standard real Gaussian random variables. He has shown that \( \mathbb E_n(\Omega) \), the expected number of zeros of \( P_n(z) \) on a measurable set \( \Omega\subset\R \), is equal to
\begin{equation}
\label{Kac}
\mathbb E_n(\Omega) = \frac1\pi\int_\Omega\frac{\sqrt{1-h_{n+1}^2(x)}}{|1-x^2|}\dd x, \quad h_{n+1}(x) = \frac{(n+1)x^n(1-x^2)}{1-x^{2n+2}},
\end{equation}
from which he proceeded with an estimate
\[
\mathbb E_n(\R) = \frac{2+o(1)}\pi\log(n+1) \quad \text{as} \quad n\to\infty.
\]
It was shown by Wilkins \cite{Wilk88}, after some intermediate results cited in \cite{Wilk88}, that there exist constants \( A_p \), \( p\geq 0 \), such that \( \mathbb E_n(\R) \) has an asymptotic expansion of the form
\begin{equation}
\label{Wilkins}
\mathbb E_n(\R) \sim \frac2\pi\log (n+1) + \sum_{p=0}^\infty A_p(n+1)^{-p}.
\end{equation}

Many recent (and already not so recent) results on random polynomials are concerned with the behavior of counting measures of zeros of random polynomials spanned by various deterministic bases with random coefficients that are not necessarily Gaussian nor i.i.d. \cite{MR1675133, MR1935565, MR2172337, MR2318650, MR2337837, MR3070481, MR3127891, MR3262481, MR3306686, LubPritXie17, MR3489554, MR3571446, MR3600134, MR3739138, MR3771497}. In the case of Kac polynomials these normalized counting measures almost surely converge to the arclength distribution on the unit circle (\(\log n \) real zeros are clearly negligible when normalized by \( 1/n \)). Our primary interest lies in studying the expected number of real zeros when the basis is a family of orthogonal polynomials in the spirit of \cite{MR0268933,MR653589,MR1377012}. More precisely, Edelman and Kostlan \cite{EdKos95} considered random functions of the form
\begin{equation}
\label{Vanderbei}
P_n(z) = \eta_0f_0(z) + \eta_1 f_1(z) + \cdots +\eta_nf_n(z),
\end{equation}
where \( \eta_i \) are certain real random variables and \( f_i(z) \) are arbitrary functions on the complex plane that are real on the real line. Using beautiful and simple geometrical argument they have shown\footnote{In fact, Edelman and Kostlan derive an expression for the real intensity function for any random vector \( (\eta_0,\ldots,\eta_n) \) in terms of its joint probability density function and of \( v(x) \).} that if \( \eta_0,\ldots,\eta_n \) are elements of a multivariate real normal distribution with mean zero and covariance matrix \( C \) and the functions \( f_i(x) \) are differentiable on the real line, then
\[
\mathbb E_n(\Omega) = \int_\Omega\rho_n(x)\dd x, \quad \rho_n(x) = \left.\frac1\pi\frac{\partial^2}{\partial s\partial t}\log\left( v(s)^\mathsf{T}Cv(t)\right)\right|_{t=s=x},
\]
where \( v(x) = \big(f_0(x),\ldots,f_n(x)\big)^\mathsf{T} \). If random variables \( \eta_i \) in \eqref{Vanderbei} are again i.i.d. standard real Gaussians, then the above expression for \( \rho_n(x) \) specializes to
\begin{equation}
\label{real-intensity}
\rho_n(x) = \frac1\pi\frac{\sqrt{K_{n+1}(x,x)K_{n+1}^{(1,1)}(x,x)-K_{n+1}^{(1,0)}(x,x)^2}}{K_{n+1}(x,x)}
\end{equation}
(this formula was also independently rederived in \cite[Proposition~1.1]{LubPritXie17} and \cite[Theorem~1.2]{uVan}), where
\[
\left\{
\begin{array}{lll}
K_{n+1}(z,w) &:=& \sum_{i=0}^nf_i(z)\overline{f_i(w)}, \medskip \\
K_{n+1}^{(1,0)}(z,w) & := & \sum_{i=0}^nf_i^\prime(z)\overline{f_i(w)}, \medskip \\
K_{n+1}^{(1,1)}(z,w) &:=& \sum_{i=0}^nf_i^\prime(z)\overline{f_i^\prime(w)}.
\end{array}
\right.
\]

In this work we concentrate on a particular subfamily of random functions \eqref{Vanderbei}, namely random polynomials of the form
\begin{equation}
\label{random-opuc}
P_n(z) = \eta_0\varphi_0(z) + \eta_1 \varphi_1(z) + \cdots +\eta_n\varphi_n(z),
\end{equation}
where \( \eta_i \) are i.i.d. standard real Gaussian random variables and \( \varphi_i(z) \) are orthonormal polynomials on the unit circle with real coefficients. That is, for some probability Borel measure \( \mu \) on the unit circle that is symmetric with respect to conjugation, it holds that
\begin{equation}
\label{opuc}
\int_\T \varphi_i(\xi)\overline{\varphi_j(\xi)}\dd\mu(\xi) = \delta_{ij},
\end{equation}
where \( \delta_{ij} \) is the usual Kronecker symbol. In this case it can be easily shown using Christoffel-Darboux formula, see \cite[Theorem~1.1]{Yatt}, that \eqref{real-intensity} can be rewritten as
\begin{equation}
\label{set-up}
\rho_n(x) = \frac1\pi\frac{\sqrt{1-h_{n+1}^2(x)}}{|1-x^2|}, \quad h_{n+1}(x) := \frac{(1-x^2)b_{n+1}^\prime(x)}{1-b_{n+1}^2(x)}, \quad b_{n+1}(x) := \frac{\varphi_{n+1}(x)}{\varphi_{n+1}^*(x)},
\end{equation}
where \( \varphi_{n+1}^*(x):=x^{n+1}\varphi_{n+1}(1/x) \) is the reciprocal polynomial (there is no need for conjugation as all the coefficients are real). When \( \mu \) is the normalized arclength measure on the unit circle, it is elementary to see that \( \varphi_m(z)=z^m \) and therefore \eqref{set-up} recovers \eqref{Kac}.

\begin{theorem}
\label{thm:1}
Let \( P_n(z) \) be given by \eqref{random-opuc}--\eqref{opuc}, where \( \mu \) is absolutely continuous with respect to the arclength measure and \( \mu^\prime(\xi) \), the respective Radon-Nikodym derivative, extends to a holomorphic non-vanishing function in some neighborhood of the unit circle. Then \( \mathbb E_n(\mu) \), the expected number of real zeros of \( P_n(z) \), satisfies
\[
\mathbb E_n(\mu) = \frac2\pi\log(n+1)+A_0 + \sum_{p=1}^{N-1}A_p^\mu(n+1)^{-p} + \mathcal O_N\left((n+1)^{-N}\right)
\]
for any integer \( N \) and all \( n \)  large, where \( \mathcal O_N(\cdot) \) depends on \( N \), but is independent of \( n \),
\[
A_0 = \frac2\pi\left(\log2+\int_0^1t^{-1}f(t)\dd t + \int_1^\infty t^{-1}(f(t)-1)\dd t\right),
\]
\( f(t) := \sqrt{1-t^2\mathrm{csch}^2t} \), and \( A_p^\mu \), \( p\geq 1 \), are some constants that do depend on \( \mu \).
\end{theorem}

Clearly, the above result generalizes \eqref{Wilkins}, where \( \dd\mu(\xi) = |\dd\xi|/(2\pi) \).

\section{Auxiliary Estimates}

In this section we gather some auxiliary estimates of quantities involving orthonormal polynomials \( \varphi_m(z) \). First of all, recall \cite[Theorem~1.5.2]{Simon1} that monic orthogonal polynomials, say \( \Phi_m(z) \), satisfy the recurrence relations
\[
\left\{
\begin{array}{l}
\displaystyle \Phi_{m+1}(z) = z\Phi_m(z) - \alpha_m\Phi_m^*(z), \medskip \\
\displaystyle \Phi_{m+1}^*(z) = \Phi_m^*(z) - \alpha_mz\Phi_m(z),
\end{array}
\right.
\]
where the recurrence coefficients \( \{\alpha_m\} \) belong to the interval \( (-1,1)  \) due to conjugate symmetry of the measure \( \mu \). In what follows we denote by \( \rho<1 \) the smallest number such that \( \mu^\prime(\xi) \) is non-vanishing and holomorphic in the annulus \( \{\rho<|z|<1/\rho\} \).

With a slight abuse of notation we shall denote various constant that depend on \( \mu \) and possibly additional parameters \( r,s \) by the same symbol \( C_{\mu,r,s} \) understanding that the actual value of \( C_{\mu,r,s} \) might be different for different occurrences, but it never depends on \( z \) or \( n \). 

\begin{lemma}
\label{lem:2}
It holds that
\[
|h_{n+1}(x)| \leq C_\mu(n+1)e^{-\sqrt{n+1}}, \quad |x| \leq 1-(n+1)^{-1/2}.
\]
\end{lemma}
\proof
Using the recurrence relations for polynomials \( \Phi_m(z) \), one can readily verify that
\[
h_{n+1}(z) = (1-z^2)\frac{(zb_n(z))^\prime}{1-(zb_n(z))^2}.
\]
Therefore, the last estimate in \cite[Section 3.3]{Yatt} implies that
\[
|h_{n+1}(x)| \leq C_\mu |\big(xb_n(x)\big)^\prime|, \quad  |x| \leq 1-(n+1)^{-1/2}.
\]
Moreover, it was shown in \cite[Equation (48)]{Yatt} that
\[
|\big(zb_n(z)\big)^\prime| \leq C_\mu(n+1)\left(r^{n-m} + \sum_{i=m}^\infty|\alpha_i| \right), \quad |z| \leq r < 1.
\]
It is further known, see \cite[Corollary~2]{MFMLS06}, that the recurrence coefficients \( \alpha_i \) satisfy
\[
|\alpha_i| \leq C_{\mu,\rho-s} s^{i+1} \quad \Rightarrow \quad \sum_{i=m}^\infty|\alpha_i| \leq \frac{C_{\mu,s-\rho} s^m}{1-\rho}, \quad \rho < s < 1,
\]
where \( C_{\mu,s-\rho} \) also depends on how close \( s \) is to \( \rho \). Given a value of the parameter \( s \), take \( m \) to be the integer part of \( -\sqrt{n+1}/\log s \) and \( r = 1 - 1/\sqrt{n+1} \). By combining the above three estimates, we deduce the desired inequality with a constant that depends on \( \mu \), \( s-\rho \), and \( s \). Optimizing the constant over \( s \) finishes the proof of the lemma.
\qed\endproof

Denote by \( D(z) \) the Szeg\H{o} function of \( \mu \), i.e.,
\[
D(z) := \exp\left\{\frac1{4\pi}\int_{\T}\frac{\xi+z}{\xi-z}\log\mu^\prime(\xi)|\dd\xi|\right\}, \quad |z|\neq1.
\]
This function is piecewise analytic and non-vanishing. Denote by \( D_{int}(z) \) the restriction of \( D(z) \) to \( |z|<1 \) and by \( D_{ext}(z) \) the restriction to \( |z|>1 \). It is known that both \( D_{int}(z) \) and \( D_{ext}(z) \) extend continuously to the unit circle and satisfy there
\[
D_{int}(\xi)/D_{ext}(\xi) = \mu^\prime(\xi), \quad |\xi|=1.
\]
Moreover, since \( \mu^\prime(\xi) \) extends to a holomorphic and non-vanishing function in the annulus \( \rho<|z|<1/\rho \), \( D_{int}(z) \) and \( D_{ext}(z) \) extend to holomorphic and non-vanishing functions in \( |z|<1/\rho \) and \( |z|>\rho \), respectively. Hence, the scattering function
\[
S(z) := D_{int}(z)D_{ext}(z), \quad \rho<|z|<1/\rho,
\]
is well defined and non-vanishing in this annulus. Since the measure \( \mu \) is conjugate symmetric, it holds that \( D(\bar z)=\overline{D(z)} \) and \( D_{ext}(1/z)=1/D_{int}(z) \). Thus, \( |S(\xi)|=1 \) for \( |\xi|=1 \) and \( S(1)=1 \). For future use let us record the following straightforward facts.

\begin{lemma}
\label{lem:3}
There exist real numbers \( s_p \), \( p\geq1 \), such that
\[
\begin{array}{rcl}
S(z) &=& 1+ \sum_{p=1}^{M-1} s_p(1-z)^p +E_M(S;z) \medskip \\
S^\prime(z) &=& - \sum_{p=0}^{M-1}(p+1)s_{p+1}(1-z)^p + E_M(S^\prime;z) \medskip \\
\log S(z) &=& \sum_{p=1}^{M-1}c_p(1-z)^p + E_M(\log S;z)
\end{array}
\]
for \( |z-1|<T<1-\rho \) and any integer \( M\geq 1 \), where the error terms satisfy
\[
\big|E_M(F;z)\big| \leq \frac{\|F\|_{|z-1|\leq T}}{1-|1-z|/T}\left(\frac{|1-z|}{T}\right)^M
\]
and \( c_p = s_p + \sum_{k=2}^p \frac{(-1)^{k-1}}k \sum_{j_1+\cdots +j_k=p}s_{j_1}\cdots s_{j_k} \). Moreover, \( s_2=s_1(s_1+1)/2 \). In particular, \( c_1=s_1 \) and \( c_2 = s_1/2 \).
\end{lemma}
\proof
Since \( c_1=s_1 \) and \( c_2 = s_2 - s_1^2/2 \), we only need to show that \( s_2=s_1(s_1+1)/2 \). It holds that \( s_1  = -S^\prime(1) \) and \( s_2 = S^{\prime\prime}(1)/2 \). Using the symmetry \( 1\equiv S(z)S(1/z) \), one can check that \( S^{\prime\prime}(1) = S^\prime(1)^2-S^\prime(1) \), from which the desired claim easily follows.
\qed\endproof

Set \( \tau:=D_{ext}(\infty) \). It has been shown in \cite[Theorem~1]{MFMLS06} that
\begin{equation}
\label{Phi-asymp}
\Phi_m(z) = \tau^{-1}z^mD_{ext}(z)\mathcal E_m(z) - \frac{\tau \mathcal I_m(z)}{D_{int}(z)}, \quad \rho<|z|<1/\rho,
\end{equation}
for some recursively defined functions \( \mathcal E_m(z), \mathcal I_m(z) \) holomorphic in the annulus \( \rho<|z|<1/\rho \) that satisfy
\begin{equation}
\label{EI-bounds}
\big|\mathcal E_m(z)-1\big| \leq \frac{C_{\mu,s}s^{2m}}{1/s-|z|} \quad \text{and} \quad \big|\mathcal I_m(z)\big| \leq \frac{C_{\mu,s}s^m}{|z|-s}, \quad \rho<s<|z|<1/s,
\end{equation}
for some explicitly defined constant \( C_{\mu,s} \), see \cite[Equations~(34)-(35)]{MFMLS06}. In particular, it follows from \eqref{Phi-asymp} that
\begin{equation}
\label{b-asymp}
b_{n+1}(z) = z^{n+1}S(z)H_n(z), \quad H_n(z) := \frac{ \mathcal{E}_{n+1}(z) - \tau^2 z^{-(n+1)} S^{-1}(z) \mathcal{I}_{n+1}(z)  }{  \mathcal{E}_{n+1}(1/z) - \tau^2 z^{n+1} S(z) \mathcal{I}_{n+1}(1/z) },
\end{equation}
for \( \rho<|z|<1/\rho \). It can be checked that the conjugate symmetry of \( \mu \) yields real-valuedness of \( H_n(z) \) on the real line. Bounds \eqref{EI-bounds} also imply that \( H_n(x) \) is close to \( 1 \) near \( x=1 \). More precisely,  the following lemma holds.

\begin{lemma}
\label{lem:4}
It holds for any \( \rho<\rho_*<1 \) that
\[
|H_n(x)-1|,|\log H_n(x)| \leq (1-x)C_{\mu,\rho_*}e^{-\sqrt{n+1}}, \quad\rho_*\leq x \leq1.
\]
Moreover, it also holds that \( |H_n^\prime(x)| \leq C_{\mu,\rho_*}e^{-\sqrt{n+1}} \) on the same interval.
\end{lemma}
\proof
Define \( W_n(z) := \mathcal{E}_{n+1}(z) - 1 - \tau^2 z^{-(n+1)} S^{-1}(z) \mathcal{I}_{n+1}(z) \) and choose \( \rho<s<s_*<\rho_*<1 \). Since \( S(z) \) is a fixed non-vanishing holomorphic function in the annulus \( \rho<|z|<1/\rho \), it follows from \eqref{EI-bounds} that
\[
|W_n(z)| \leq C_{\mu,s,s_*}\big(s/s_*\big)^n, \quad s_*\leq |z|\leq1/s_*.
\]
It further follows from the maximum modulus principle that
\[
|W_n(z)-W_n(1/z)| \leq |1-z|C_{\mu,s,s_*}\big(s/s_*\big)^n, \quad s_*\leq |z|\leq1/s_*,
\]
where, as agreed before, the actual constants in the last two inequalities are not necessarily the same. Since \(|\log(1+\zeta)|\leq 2|\zeta| \) for \( |\zeta|\leq 1/2 \), there exists a constant  \( A_{\mu,s,s_*} \) such that
\[
|H_n(z)-1|,| \log H_n(z)| \leq |1-z|A_{\mu,s,s_*}\big(s/s_*\big)^n, \quad s_*\leq |z|\leq1/s_*.
\]
Observe that the constants \( A_{\mu,s,s_*}e^{\sqrt{n+1}}\big(s/s_*\big)^n \) are uniformly bounded above. Then the first claim of the lemma follows by minimizing these constants over all parameters \( s<s_* \) between \( \rho \) and \( \rho_* \). Further, it follows from Cauchy's formula that
\[
H_n^\prime(z) = \left(\int_{|\zeta|=1/s_*}-\int_{|\zeta|=s_*}\right)\frac{H_n(\zeta)-1}{(\zeta-z)^2}\frac{\dd \zeta}{2\pi\ic}
\] 
for \( \rho_* \leq |z| \leq 1/\rho_* \) and therefore it holds in this annulus that
\[
|H_n^\prime(z)| \leq C_{\mu,s,s_*,\rho_*}\big(s/s_*\big)^n.
\]
The last claim of the lemma is now deduced in the same manner as the first one.
\qed\endproof

\section{Proof of Theorem~\ref{thm:1}}
\label{sec:3}

Using \eqref{set-up}, it is easy to show that
\[
\mathbb E_n(\mu) = \frac2\pi\int_{-1}^1\frac{\sqrt{1-h_{n+1}^2(x)}}{1-x^2}\dd x.
\]
Furthermore, if we define \( \dd\sigma(\xi) := \mu^\prime(-\xi)|\dd \xi| \), then \( \sigma^\prime(\xi) = \mu^\prime(-\xi) \) is still holomorphic and positive on the unit circle. Moreover, \( b_n(z;\sigma) = b_n(-z;\mu) \). Therefore,
\begin{equation}
\label{reduction}
\mathbb E_n(\mu) = \widehat{\mathbb E}_n(\mu) +  \widehat{\mathbb E}_n(\sigma), \quad  \widehat{\mathbb E}_n(\nu) := \frac2\pi\int_0^1\frac{\sqrt{1-h_{n+1}^2(x;\nu)}}{1-x^2}\dd x,
\end{equation}
for \( \nu\in\{\mu,\sigma\} \). Thus, it is enough to investigate the asymptotic behavior of \(  \widehat{\mathbb E}_n(\mu) \). To this end, let
\begin{equation}
\label{notation}
a:=(n+1)^{1/2} \quad \text{and} \quad x =: 1-t/(n+1), \quad 0\leq t\leq a.
\end{equation}
We shall also write
\begin{equation}
\label{En}
1 - h_{n+1}^2(x) =: f^2(t)(1+E_n(t)),
\end{equation}
for \( 1-(n+1)^{-1/2}\leq x\leq 1 \), where \( f(t) \) was defined in Theorem~\ref{thm:1}.

\begin{lemma}
\label{lem:5}
Given an integer \( N\geq 1 \), it holds that
\[
\widehat{\mathbb E}_n(\mu) = \frac1\pi\log(n+1) + \frac12A_0 + G_n(t) - \frac12\sum_{p=1}^{N-1}H_p(n+1)^{-p} + \mathcal O_N\left((n+1)^{-N}\right)
\]
for large \( n \), where \( \mathcal O_N(\cdot) \) is independent of \( n \), but does depend on \( N \),
\[
G_n(t) := \frac1\pi\int_0^a\left(t^{-1} + \big(2(n+1)-t\big)^{-1}\right)f(t)\left(\big(1+E_n(t)\big)^{1/2}-1\right)\dd t,
\]
and \( H_p := \displaystyle  \frac1{2^{p-1}\pi}\int_0^\infty\big( 1-f(t) \big)t^{p-1}\dd t \) for \( p\geq1 \).
\end{lemma}
\proof
Set \( \delta:=1-(n+1)^{-1/2} \). It trivially holds that
\[
\widehat{\mathbb E}_n(\mu) = \frac{2}{\pi} \int_{0}^{\delta} \frac{\dd x}{1-x^2} - \frac2\pi\int_0^\delta\frac{1-\sqrt{1-h_{n+1}^2(x)}}{1-x^2}\dd x + \frac2\pi\int_\delta^1\frac{\sqrt{1-h_{n+1}^2(x)}}{1-x^2}\dd x.
\]
Denote the third integral above by \( B_n(t) \). The second integral above is positive and equals to
\[
\frac2\pi\int_0^\delta\frac{h_{n+1}^2(x)}{1+\sqrt{1-h_{n+1}^2(x)}}\frac{\dd x}{1-x^2} \leq \frac2\pi\int_0^\delta h_{n+1}^2(x)\frac{\dd x}{1-\delta^2} =\mathcal O\left(a^5e^{-2a}\right),
\]
where we used Lemma~\ref{lem:2} for the last estimate. Therefore,
\[
\widehat{\mathbb E}_n(\mu)  = \frac{1}{\pi} \log\left(\frac{1+\delta}{1-\delta}\right) + B_n(t) + o_N\left((n+1)^{-N}\right),
\]
where \( o_N(\cdot) \) is independent of \( n \), but does depend on \( N \). Substituting \( x = 1 - t/(n+1) \) into the expression for \( B_n(t) \) and recalling \eqref{En}, we get that
\begin{eqnarray*}
B_n(t) & = & \frac{1}{\pi} {\displaystyle \int_{0}^{a} f(t)\big( 1+E_n(t) \big)^{1/2} \frac{2(n+1)}{t(2(n+1)-t)}\,\dd t} \\
& = & \frac{1}{\pi} \left( \log 2 + \log \frac{1}{1 + \delta} \right) + \frac{1}{\pi} \int_{0}^{a} \frac{f(t)}{t}\,\dd t - \frac{1}{\pi}  \int_{0}^{a}\frac{1-f(t)}{2(n+1)-t} \,\dd t + G_n(t).
\end{eqnarray*}
It was shown in \cite[Lemma~8]{Wilk88} that
\[
\frac{1}{\pi} {\displaystyle \int_{0}^{a}\frac{1-f(t)}{2(n+1)-t} \,\dd t} = \frac12\sum_{p=1}^{N-1}H_p(n+1)^{-p} + \mathcal O_N\left((n+1)^{-N}\right),
\]
where \( \mathcal O_N(\cdot) \) is independent of \( n \), but does depend on \( N \). Moreover, it holds that
\begin{multline*}
\frac{1}{\pi} \log\left(\frac{1+\delta}{1-\delta}\right) + \frac{1}{\pi} \left( \log 2 + \log \frac{1}{1 + \delta} \right) + \frac{1}{\pi} {\displaystyle \int_{0}^{a} \frac{f(t)}{t}\,\dd t} = \\ = \frac{1}{\pi} \log \frac{a}{1-\delta} + \frac{1}{2} A_{0} + \frac{1}{\pi} {\displaystyle \int_{a}^{\infty} \frac{1-f(t)}{t}\,\dd t} .
\end{multline*}
Since \(  \log a - \log(1-\delta) = \log (n+1) \) and it was shown in \cite[Lemma~7]{Wilk88} that
\[
\frac{1}{\pi} {\displaystyle \int_{a}^{\infty} \frac{1-f(t)}{t}\,\dd t}  = \mathcal O\left(a e^{-2a}\right)= o_N\left( (n+1)^{-N} \right),
\]
where as usual \( o_N(\cdot) \) is independent of \( n \), but does depend on \( N \), the claim of the lemma follows.
\qed\endproof

We continue by deriving a different representation for the functions \( E_n(t) \). To this end, notice that \( t^2\mathrm{csch}^2t = 1 -t^2/3 + \mathcal O\big(t^4\big) \) as \( t\to 0 \) and therefore \( f^2(t) = t^2/3 + \mathcal O\big(t^4\big) \) as \( t\to 0 \). Hence, the function
\begin{equation}
\label{chi}
\chi(t):= \left(\frac{t^2\mathrm{csch}t}{f(t)}\right)^2
\end{equation}
is continuous and non-vanishing at zero. Once again, we use notation from \eqref{notation}.

\begin{lemma}
\label{lem:6}
Set \( b_{n+1}^2(x) =: e^{-\mu_n(t)-2t} \) and \( b_{n+1}^\prime(x) =: (n+1)e^{w_n(t)-t} \). Then it holds that
\[
E_n(t) = t^{-2}\chi(t)\left[1-\left(1-\frac t{2(n+1)}\right)^2\frac{e^{2w_n(t)}}{(1+D_n(t))^2}\right], \quad D_n(t) := \frac{1-e^{-\mu_n(t)}}{e^{2t}-1}.
\]
Moreover, \( \lim_{t\to0^+}E_n(t) \) exists and is finite.
\end{lemma}
\proof
Since \( h_{n+1}(1)=1 \) and \( x = 1-t/(n+1) \), it follows from \eqref{En} and the L'H\^opital's rule that
\[
\lim_{t\to0^+}E_n(t) = \frac6{(n+1)^2}\lim_{x\to1^-}\frac{1-h_{n+1}(x)}{(1-x)^2} - 1 = \frac3{(n+1)^2}\lim_{x\to1^-}\frac{h_{n+1}^\prime(x)}{1-x} - 1.
\]
Since \( h_{n+1}(z) \) is a holomorphic function around \( 1 \), the latter limit is finite if and only if \( h_{n+1}^\prime(1) =0 \). As Blaschke products \( b_{n+1}(z) \) satisfy \( b_{n+1}(x)b_{n+1}(1/x)\equiv1 \), it holds that \( h_{n+1}(x)=h_{n+1}(1/x) \), which immediately yields the desired equality.

To derive the claimed representation of \( E_n(t) \), recall \eqref{set-up} and substitute \( x = 1-t/(n+1) \) into \eqref{En} to get that
\[
\begin{split}
f^2(t)(1+E_n(t))  & = 1 - \left( 1-\frac{t}{2(n+1)} \right)^2 \frac{4 t^2 e^{2w_n(t)-2t}}{ \big(1-e^{-\mu_n(t)-2t}\big)^2}\\
& = 1 - \left( 1-\frac{t}{2(n+1)} \right)^2 \frac{t^2 \mathrm{csch}^2t  e^{2w_n(t)}}{ \big(1+D_n(t)\big)^2}\\
& = f^2(t) \left[ 1 + t^{-2}\chi(t)  \left(1 - \left( 1-\frac{t}{2(n+1)} \right)^2 \frac{e^{2w_n(t)}}{ \big(1+D_n(t)\big)^2} \right)\right]
\end{split} 
\]
from which the first claim of the lemma easily follows.
\qed\endproof
 
In the next four lemmas we repeatedly use approximation by Taylor polynomials with the Lagrange remainder:
\begin{equation}
\label{taylor}
F(y) = \sum_{k=0}^{M-1}\frac{F^{(k)}(0)}{k!}y^K + \frac{F^{(M)}(\theta y)}{M!}y^M
\end{equation}
for some \( \theta\in(0,1) \) that dependents on both \( y \) and \( M \).
 
\begin{lemma}
\label{lem:7}
Put \( \omega(t) := t/(e^{2t}-1) \). Given an integer \( N\geq 1 \), it holds for all \( n \) large that
\[
\big(1+D_n(t)\big)^{-2}  = 1 + \sum_{p=1}^{N-1}\alpha_p(t)(n+1)^{-p} + \alpha_{n,N}(t)(n+1)^{-N},
\]
where the functions \( \alpha_p(t) \) are independent of $n$ and $N$ and are polynomials of degree \( p \) in \( \omega \) with coefficients that are polynomials in \( t \) of degree at most \( 2p-1 \), and the functions \( \alpha_{n,N}(t) \) are bounded in absolute value for \( 0\leq t\leq a \) by a polynomial of degree $ 2N-1 $ whose coefficients are independent of $n$. Moreover,
\[
\alpha_p(t) = (p+1)s_1^p - ps_1^{p-1}(2s_1+1)t + \mathcal O\big(t^2\big) \quad \text{as} \quad t\to 0.
\]
\end{lemma}
\proof
We start by deriving an asymptotic expansion of \( \mu_n(t) \). It follows from Lemma~\ref{lem:4} that \( \log H_n(x)=t\mathcal O(a^{-2}e^{-a}) = to_N(1)(n+1)^{-N} \) uniformly for \( 0\leq t\leq a \). Fix \( T \) in Lemma~\ref{lem:3} and let \( n_T \) be such that \( 1<\sqrt{n_T+1}T \). Then it holds for all \( n\geq n_T \) that
\[
\log (SH_n)(x) = \sum_{p=1}^{N-1}c_p t^p(n+1)^{-p} + t\hat c_N(t)(n+1)^{-N},
\]
where \( |\hat c_N(t)| \leq C_{\mu,T,N}t^{N-1} + o_N(1)\) uniformly for \( 0\leq t\leq a \) and \( C_{\mu,T,N}\leq C_{\mu,T} T^{-N} \). Hence, it follows from \eqref{b-asymp} and \cite[Lemma~2]{Wilk88} that
\begin{eqnarray}
\mu_n(t) & = & -2(n+1)\log x - 2t - 2\log(SH_n)(x) \nonumber \\
\label{mun1}
& = & \sum_{p=1}^{N-1}t^pm_p(t)(n+1)^{-p} + tm_{n,N}(t)(n+1)^{-N},
\end{eqnarray}
where
\[
m_p(t):=\big(2(p+1)^{-1}t-2c_p\big) \quad \text{and} \quad m_{n,N}(t) := 2 \hat m_{n,N}(t) t^N/(N+1)-2\hat c_N(t)
\]
with \( 1 \leq \hat m_{n,N}(t) \leq (3/2)^{N+1} \). Assuming that \( T<2/3 \), we have that
\begin{equation}
\label{mun3}
|m_{n,N}(t)| \leq C_{\mu,T,N}t^{N-1}(t+1) + o_N(1)
\end{equation}
uniformly for \( 0\leq t\leq a \) and \( C_{\mu,T,N}\leq C_{\mu,T} T^{-N} \). Using \eqref{mun1} with \( N=1 \), we get that
\begin{equation}
\label{mun2}
|\mu_n(t)| = \left|\frac{tm_{n,1}(t)}{n+1}\right| \leq \frac{|m_{n,1}(t)|}{\sqrt{n+1}} \leq C_{\mu,T}, \quad 0\leq t\leq a.
\end{equation}

Recalling the definition of \( D_n(t) \) in Lemma~\ref{lem:6}, we get from \eqref{taylor} that 
\[
D_n(t) = \omega(t)\frac{1-e^{-\mu_n(t)}}t = \omega(t)\left(-\frac1t\sum_{k=1}^{N-1}\frac{(-1)^k}{k!}\mu_n^k(t) -\frac1t e^{-\theta_1\mu_n(t)}\frac{(-1)^N}{N!}\mu_n^N(t)\right)
\]
for some \( \theta_1\in(0,1) \) that depends on \( N \) and \( \mu_n(t) \). Plugging \eqref{mun1} into the above formula gives us
\begin{equation}
\label{dn1}
D_n(t) = \omega(t)\sum_{p=1}^{N-1}t^{p-1}d_p(t)(n+1)^{-p} + \omega(t)d_{n,N}(t)(n+1)^{-N},
\end{equation}
where \( d_p(t) \) is a polynomial of degree \( p \) with coefficients independent of \( n \) and \( N \) given by
\[
d_p(t) := -\sum_{k=1}^p\frac{(-1)^k}{k!}\sum_{j_1+\cdots+j_k=p}m_{j_1}(t)\cdots m_{j_k}(t),
\]
here, each index \( j_i\in\{1,\ldots,p\} \), and \( d_{n,N}(t) \) is given by
\[
d_{n,N}(t) := -\sum_{k=1}^{N-1}\frac{(-1)^k}{k!}\sum_{j_1+\cdots+j_k\geq N}\frac1t\frac{m_{n,j_1,N}(t)\cdots m_{n,j_k,N}(t)}{(n+1)^{j_1+\cdots+j_k-N}} - \frac{(-1)^N}{N!}\frac{(n+1)^N}{e^{\theta_1\mu_n(t)}}\frac{\mu_n^N(t)}t
\]
with \( m_{n,j,N}(t) := t^jm_j(t) \) when \( j<N \) and \( m_{n,N,N}(t) := tm_{n,N}(t) \). Recall that \( t^2/(n+1)\leq 1 \) on \( 0\leq t\leq a \) since \( a=\sqrt{n+1} \). Hence, the first summand above is bounded in absolute value for \( 0\leq t\leq a \) by a polynomial of degree \( 2N-1 \) whose coefficients depend on \( N \) but are independent of \( n \). We also get from \eqref{mun2} and \eqref{mun3} that
\[
\left|e^{-\theta_1\mu_n(t)}(n+1)^N\mu_n^N(t)/t\right| \leq e^{C_{\mu,T}}t^{N-1}|m_{n,1}(t)|^N \leq C_{\mu,T}^*t^{N-1}(t+2)^N
\]
for \( 0\leq t\leq a \). Further, using \eqref{dn1} with \( N=1 \) and \eqref{mun2} gives us
\begin{equation}
\label{dn2}
|D_n(t)| = \frac{\omega(t)}{e^{\theta_1\mu_n(t)}}\left|\frac{\mu_n(t)}t\right| \leq \frac{e^{C_{\mu,T}}}2\frac{|m_{n,1}(t)|}{n+1} \leq \frac{C_{\mu,T}e^{C_{\mu,T}}}{2\sqrt{n+1}}, \quad 0\leq t\leq a.
\end{equation}
Notice also that since \( c_1=s_1 \) and \( c_2= s_1/2 \) by Lemma~\ref{lem:3}, we have that
\[
d_1(t) = t -2s_1 \quad \text{and} \quad d_2(t) = -(1/2)t^2 + t(2s_1+2/3) - s_1(2s_1+1).
\]

It follows from \eqref{dn2} that for any \( -1<D<0 \), there exists an integer \( n_D\geq n_T \) such that \( D\leq D_n(t) \) for \( 0\leq t\leq a \) and \( n \geq n_D \). Hence, we get from \eqref{taylor} that
\[
\big(1+D_n(t)\big)^{-2}  = 1 + \sum_{k=1}^{N-1} (-1)^k(k+1) D_n^k(t) + \frac{(-1)^N(N+1)D_n^N(t)}{(1+\theta_2D_n(t))^{N+2}}
\]
for all \( n\geq n_D \) and some \( \theta_2\in(0,1) \) that depends on \( N \) and \( D_n(t) \). Then the statement of the lemma follows with
\[
\alpha_p(t) := \sum_{k=1}^p(-1)^k(k+1)\omega^k(t)t^{p-k}\sum_{j_1+\cdots+j_k=p}d_{j_1}(t)\cdots d_{j_k}(t)
\]
here again, each index \( j_i\in\{1,\ldots,p\} \),  and
\[
\alpha_{n,N}(t) := \sum_{k=1}^{N-1}(-1)^k(k+1)\omega^k(t)\sum_{j_1+\cdots+j_k\geq N}\frac{d_{n,j_1,N}(t)\cdots d_{n,j_k,N}(t)}{(n+1)^{j_1+\cdots+j_k-N}} + (n+1)^N\frac{(-1)^N(N+1)D_n^N(t)}{(1+\theta_2D_n(t))^{N+2}}
\]
with \( d_{n,j,N}(t) := t^{j-1}d_j(t) \) when \( j<N \) and \( d_{n,N,N}(t) := d_{n,N}(t) \). Reasoning as before lets us conclude that the first summand in the definition of \( \alpha_{n,N}(t) \) is bounded in absolute value for \( 0 \leq t \leq a \) by a polynomial of degree \( 2N-1 \) whose coefficients depend on \( N \) but are independent of \( n \). Moreover, since
\[
\left|\frac{(n+1)^ND_n^N(t)}{(1+\theta_2D_n(t))^{N+2}}\right| \leq \frac {e^{NC_{\mu,T}}|m_{n,1}(t)|^N}{2^N(1-D)^{N+2}} \leq \frac {C_{\mu,T}^*e^{NC_{\mu,T}}(t+2)^N}{2^N(1-D)^{N+2}}, \quad 0\leq t\leq a,
\]
by \eqref{dn2} and \eqref{mun3}, the same is true for the second summand as well. Now, notice that
\[
\alpha_p(t) = \big(-\omega(t)d_1(t) \big)^{p-2}\left( (p+1)\big(\omega(t)d_1(t)\big)^2 -p(p-1)t\omega(t) d_2(t)\right) + \mathcal O\big(t^2\big)
\]
as \( t\to 0 \). Since \( 2\omega(t) = 1-t + \mathcal O\big(t^2\big) \) as \( t\to 0 \), the last claim of the lemma follows after a  straightforward computation.
\qed\endproof

\begin{lemma}
\label{lem:8}
Given \( N\geq 1 \), it holds for all \( n \) large that
\[
e^{2w_n(t)} = 1 + \sum_{p=1}^{N-1} \beta_p(t)(n+1)^{-p} + \beta_{n,N}(t)(n+1)^{-N},
\]
where \( \beta_p(t) \) is a polynomial of degree \( 2p \) whose coefficients are independent of $n$ and $N$ and the functions \( \beta_{n,N}(t) \) are bounded in absolute value when \( 0 \leq t \leq a \) by a polynomial of degree $ 2N$ whose coefficients are independent of $n$. Moreover, as \( t\to0 \), it holds that
\[
\begin{cases}
\beta_1(t) = -2s_1 + 2(s_1+1)t - t^2, \\
\beta_2(t) = s_1^2 - 4s_1(s_1+1)t + \mathcal O\big(t^2\big),\\
\beta_3(t) = 2s_1^2(s_1+1)t+ \mathcal O\big(t^2\big),\\
\beta_p(t) = \mathcal O\big(t^2\big), \quad p\geq 4.
\end{cases}
\]
\end{lemma}
\proof
We start by deriving an asymptotic expansion for \( w_n(t) \).  It follows from the very definition of \( w_n(t) \) in Lemma~\ref{lem:6}, \eqref{b-asymp}, and \cite[Lemma~2]{Wilk88} that
\begin{eqnarray*}
w_n(t) & = & t + \log \frac{b_{n+1}'(x)}{n+1} = t + n \log x + \log\left( (SH_n)(x)+\frac{x(SH_n)'(x)}{n+1} \right) \\
& = & \sum_{p=1}^{N-1} t^p\phi_{p}(t) (n+1)^{-p} + \phi_{n,N}(t) (n+1)^{-N} + \log\left( (SH_n)(x)+\frac{x(SH_n)'(x)}{n+1} \right),
\end{eqnarray*}
where
\begin{equation}
\label{phip}
\phi_{p}(t):=\frac{p+1-pt}{p(p+1)} \quad \text{and} \quad \phi_{n,N}(t) := \left( N^{-1} - \frac{n \hat m_{n,N}(t)t}{(N+1)(n+1)}\right)t^N
\end{equation}
with some \( 1\leq \hat m_{n,N}(t) \leq(3/2)^N \). Further, notice that
\[
(S^{(i)}H_n)(x) = S^{(i)}(x) + o_N(1)(n+1)^{-N} \quad \text{and} \quad (SH_n^\prime)(x) = o_N(1)(n+1)^{-N}
\]
uniformly for \( 0\leq t\leq a \), \( i\in\{0,1\} \), by Lemma~\ref{lem:4} and since \( S(z) \) is a fixed holomorphic function in a neighborhood of \( 1 \). Fix \( T \) in Lemma~\ref{lem:3}. Then it holds for all \( n\geq n_T \) that
\[
(SH_n)(x) = 1 + \sum_{j=1}^{N-1} s_j \frac{t^j}{(n+1)^j} + \hat s_N(t) (n+1)^{-N},
\]
and
\[
(SH_n)'(x) = - \sum_{j=1}^{N-1} j s_j \frac{t^{j-1}}{(n+1)^{j-1}} - \hat f_N(t) (n+1)^{-N},
\]
where  \( |\hat s_N(t)|,|\hat f_N(t)| \leq C_\mu(t/T)^N + o_N(1) \) uniformly for \( 0\leq t\leq a \). Therefore,
\begin{equation}
\label{ln1}
L_n(t) := (SH_n)(x)-1+\frac{x(SH_n)'(x)}{n+1} = \sum_{j=1}^{N-1} t^{j-1} l_j(t)(n+1)^{-j}  + l_{n,N}(t) (n+1)^{-N},
\end{equation}
where
\[
l_j(t) := \big( s_j (t - j) + (j-1) s_{j-1} \big)
\]
and
\[
l_{n,N}(t) := (N-1)s_{N-1} t^{N-1} + \hat s_N(t) - \left(1-\frac{t}{n+1}\right)\frac{\hat f_N(t)}{n+1}.
\]
In particular, it holds that
\begin{equation}
\label{ln3}
|l_{n,N}(t)| \leq 2C_\mu(t/T)^N + (N-1)s_{N-1} t^{N-1} + o_N(1)
\end{equation}
and therefore
\begin{equation}
\label{ln2}
|L_n(t)| \leq \frac{|l_{n,1}(t)|}{n+1} \leq \frac{C_{\mu,T}}{\sqrt{n+1}}, \quad 0\leq t\leq a.
\end{equation}
Hence, given \( -1<L<0 \), there exists an integer \( n_L\geq n_T \) such that \( L\leq L_n(t) \) for \( 0\leq t\leq a \) and \( n \geq n_L \). Thus, we get from \eqref{taylor} that
\[
\log(1+L_n(t)) = \sum_{k=1}^{N-1}\frac{(-1)^{k-1}}{k}L_n^k(t) + \frac{(-1)^{N-1}L_n^N(t)}{N(1+\theta_3 L_n(t))^N}
\]
for some \( \theta_3\in(0,1) \) that depends on \( N \) and \( L_n(t) \). Therefore, we get from \eqref{ln1} that
\[
\log\left( (SH_n)(x)+\frac{x(SH_n)'(x)}{n+1} \right) = \sum_{p=1}^{N-1} \psi_{p}(t) (n+1)^{-p} + \psi_{n,N}(t) (n+1)^{-N},
\]
where \( \psi_p(t) \) is a polynomial of degree \( p \) with coefficients independent of \( n \) and \( N \) given by
\begin{equation}
\label{psip}
\psi_p(t) := \sum_{k=1}^p\frac{(-1)^{k-1}}{k}\sum_{j_1+\cdots+j_k=p} t^{p-k} l_{j_1}(t)\cdots l_{j_k}(t),
\end{equation}
here, each index \( j_i\in\{1,\ldots,p\} \), and \( \psi_{n,N}(t) \) is given by
\[
\psi_{n,N}(t) := \sum_{k=1}^{N-1}\frac{(-1)^{k-1}}{k}\sum_{j_1+\cdots+j_k\geq N} \frac{l_{n,j_1,N}(t)\cdots l_{n,j_k,N}(t)}{(n+1)^{j_1+\cdots+j_k-N}} + (n+1)^N\frac{(-1)^{N-1}L_n^N(t)}{N(1+\theta_3 L_n(t))^N}
\]
with \( l_{n,j,N}(t):= t^{j-1} l_j(t) \) when \( j<N \) and \( l_{n,N,N}(t):= l_{n,N}(t) \). As in the previous lemma, since \( t^2/(n+1)\leq 1 \) when \( 0\leq t\leq a \), the first summand above is bounded in absolute value by a polynomial of degree \( N \) whose coefficients are independent of \( n \). It also follows from \eqref{ln2} and \eqref{ln3} that
\[
\frac{(n+1)^N|L_n^N(t)|}{|1+\theta_3 L_n(t)|^N} \leq \frac{|l_{n,1}(t)|^N}{(1-L)^N} \leq C_{\mu,T}\frac{(t+1)^N}{(1-L)^N}, \quad 0\leq t\leq a,
\]
for all \( n\geq n_L\). Altogether, we have shown that
\begin{equation}
\label{wn1}
w_n(t) = \sum_{p=1}^{N-1} \big( t^p \phi_{p}(t) + \psi_p(t) \big) (n+1)^{-p} + \big( \phi_{n,N}(t) + \psi_{n,N}(t) \big) (n+1)^{-N} 
\end{equation}
with \( \phi_p,\psi_p \) and \( \phi_{n,N},\psi_{n,N} \) as described above. We also can deduce from \eqref{phip} and \eqref{psip} that \( t\phi_1(t) + \psi_1(t) = -s_1 + t (s_1+1) - t^2/2 \) and
\begin{equation}
\label{psip-l}
t^p\phi_p(t) + \psi_p(t)  = \frac{(-1)^{p-1}}{p} l_{1}^p(t) + (-1)^{p-2} t l_{1}^{p-2}(t)l_{2}(t) + \mathcal{O}\big(t^2\big) = - \frac{s_1^p}{p} + \mathcal{O}\big(t^2\big)
\end{equation}
for \( p \geq 2 \), where we used that \( 2s_2 = s_1^2+ s_1 \), see Lemma~\ref{lem:3}. Since
\[
\big|\psi_{n,1}(t) \big| \leq  (n+1)\frac{|L_n(t)|}{1-L} \leq \sqrt{n+1}\frac{C_{\mu,T}}{1-L}, \quad 0\leq t\leq a,
\]
by \eqref{ln2} for \( n\geq n_L \), we get from \eqref{wn1}, applied with \( N=1 \), and \eqref{phip} that
\begin{equation}
\label{wn2}
|w_n(t)| = \left|\frac{\phi_{n,1}(t) + \psi_{n,1}(t)}{n+1}\right| \leq C_{\mu,T,L}, \quad 0\leq t\leq a, \quad n\geq n_L.
\end{equation}

Now, using \eqref{taylor} once more, we get
\[
e^{2w_n(t)} = 1 + \sum_{k=1}^{N-1}\frac{2^k}{k!}w_n^k(t) + e^{2\theta_4w_n(t)}\frac{(2)^N}{N!}w_n^N(t)
\]
for some \( \theta_4\in(0,1) \) that depends on \( N \) and \( w_n(t) \). Plugging \eqref{wn1} into the above formula gives us the desired expansion with
\begin{equation}
\label{ap}
\beta_p(t) := \sum_{k=1}^p\frac{2^k}{k!}\sum_{j_1+\cdots+j_k=p} \big( t^{j_1}\phi_{j_1}(t) + \psi_{j_1}(t) \big)\cdots 
\big( t^{j_k} \phi_{j_k}(t) + \psi_{j_k}(t) \big),
\end{equation}
which is a polynomial of degree \( 2p \) with coefficients independent of \( n \) and \( N \), and
\[
\beta_{n,N}(t) := \sum_{k=1}^{N-1}\frac{2^k}{k!}\sum_{j_1+\cdots+j_k\geq N} \frac{ \prod_{i=1}^k\big( \phi_{n,j_i,N}(t) + \psi_{n,j_i,N}(t)\big)}{(n+1)^{j_1+\cdots+j_k-N}} + e^{2\theta_4w_n(t)}\frac{2^N}{N!}(n+1)^N w_n^N(t)
\]
with \( \phi_{n,j,N}(t) := t^j \phi_{j}(t), \psi_{n,j,N}(t) := \psi_{j}(t) \) when \( j<N \) and \( \phi_{n,N,N}(t) := \phi_{n,N}(t), \psi_{n,N,N}(t) := \psi_{n,N}(t) \), which is bounded in absolute value when \( 0\leq t\leq a \)  by a polynomial of degree \( 2N \) whose coefficients are independent of \( n \) due to \eqref{wn2} and the same reasons as in the similar previous computations. Thus, it only remains to compute the linear approximation to \( \beta_p(t) \) at zero. Now, it follows from \eqref{psip-l} and \eqref{ap} that
\begin{multline*}
\beta_p(t) = s_1^p \sum_{k=1}^p \frac{(-2)^k}{k!} \sum_{j_1+\cdots+j_k=p} \frac{1}{j_1 \cdots j_k}  \\ -  \left(s_1^{p-1}(s_1+1) \sum_{k=1}^p \frac{(-2)^k}{k!} \sum_{j_1+\cdots+j_k=p} \frac{n(j_1,\ldots,j_k)}{j_1 \cdots j_k}\right) t + \mathcal O\big(t^2\big) 
\end{multline*}
where \( n(j_1,\ldots,j_k) \) is the number of \( 1 \)'s in the partition \( \{j_1,\ldots,j_k\} \) of \( p \). To simplify this expression observe that
\begin{equation}
\label{trick}
\begin{split}
(1-x)^2e^{-2yx} & = e^{2\log{(1-x)}-2yx} = 1 + \sum_{k=1}^{\infty} \frac{(-2)^k}{k!}\big(yx - \ln{(1-x)} \big)^k\\
& =  1 + \sum_{k=1}^{\infty} \frac{(-2)^k}{k!}\left((1+y)x + \sum_{j=2}^\infty \frac{x^j}j \right)^k \\
& = 1 + \sum_{p=1}^{\infty} \left( \sum_{k=1}^p \frac{(-2)^k}{k!} \sum_{j_1+\cdots+j_k=p} \frac{(1+y)^{n(j_1,\ldots,j_k)}}{j_1 \cdots j_k} \right) x^p,
\end{split}
\end{equation}
where \( y \) is a free parameter. By putting \( y=0 \) in this expression, we get that
\[
\sum_{k=1}^p \frac{(-2)^k}{k!} \sum_{j_1+\cdots+j_k=p} \frac{1}{j_1 \cdots j_k} = 
\left\{
\begin{array}{rl}
-2 & \text{if} \quad  p = 1, \\
1  & \text{if} \quad  p = 2, \\
0 & \text{if} \quad  p \geq 3.
\end{array}
\right.
\]
Moreover, by differentiating \eqref{trick} with respect to \( y \) and then putting \( y=0 \), we get
\[
\sum_{k=1}^p \frac{(-2)^k}{k!} \sum_{j_1+\cdots+j_k=p} \frac{n(j_1,\ldots,j_k)}{j_1 \cdots j_k} = 
\left\{
\begin{array}{rl}
-2 & \text{if} \quad  p = 1, \\
4  & \text{if} \quad  p = 2, \\
-2 & \text{if} \quad  p = 3, \\
0 &  \text{if} \quad  p \geq 4,
\end{array}
\right.
\]
which clearly finishes the proof of the last claim of the lemma.
\qed\endproof

\begin{lemma}
\label{lem:9}
Let \( \chi(t) \) be given by \eqref{chi}. For any integer \( N\geq 1 \), it holds that
\[
\big(1+E_n(t)\big)^{1/2}-1 = \chi(t)\sum_{p=1}^{N-1} u_p(t)(n+1)^{-p} + \chi(t)u_{n,N}(t)(n+1)^{-N},
\]
where \( u_p(t) \) is bounded in absolute value\footnote{In fact, \( u_p(t) \) is a multivariate polynomial in \( \omega,\chi \), and \( t \).} on \( 0 \leq t< \infty \) by a polynomial of degree \( 2p-2 \) whose coefficients are independent of $n$ and $N$ and the functions \( u_{n,N}(t) \) are bounded in absolute value when \( 0 \leq t \leq a \) by a polynomial of degree $ 2N-2$ whose coefficients are independent of $n$.
\end{lemma}
\proof
Set
\[
R_n(t) := \left( 1- \frac{t}{2(n+1)}\right)^2 \frac{e^{2w_n(t)}}{(1+D_n(t))^2}.
\]
Lemmas~\ref{lem:7} and~\ref{lem:8} yield that \( R_n(t) \) has the following asymptotic expansion:
\[
R_n(t) = 1 + \sum_{p=1}^{N-1} r_{p}(t)(n+1)^{-p} + r_{n,N}(t)(n+1)^{-N},
\]
where
\[
r_p(t) := \sum_{j=0}^p \beta_j(t)\alpha_{p-j}(t) - \sum_{j=0}^{p-1} t\beta_j(t)\alpha_{p-1-j}(t) + \sum_{j=0}^{p-2} t^2\beta_j(t)\alpha_{p-2-j}(t)/4
\]
with \( \alpha_0(t)=\beta_0(t) :\equiv 1 \), and \( r_{n,N}(t) \) given by
\[
\sum_{k=N}^{2N+2}\left(\sum_{j=0}^k \frac{\beta_{n,j,N}(t)\alpha_{n,k-j,N}(t)}{(n+1)^{k-N}} - \sum_{j=0}^{k-1}  \frac{t\beta_{n,j,N}(t)\alpha_{n,k-1-j,N}(t)}{(n+1)^{k-N}} + \sum_{j=0}^{k-2}  \frac{t^2\beta_{n,j,N}(t)\alpha_{n,k-2-j,N}(t)/4}{(n+1)^{k-N}} \right)
\]
with \( \alpha_{n,j,N}(t) := \alpha_{j}(t),\beta_{n,j,N}(t) := \beta_{j}(t)  \) when \( j < N \), \( \alpha_{n,N,N}(t) := \alpha_{n,N}(t), \beta_{n,N,N}(t) := \beta_{n,N}(t)  \), and \( \alpha_{n,j,N}(t) = \beta_{n,j,N}(t) :\equiv 0 \) when \( j>N \). It also follows from Lemmas~\ref{lem:7} and~\ref{lem:8} that the functions \( r_p(t) \) are independent of \(n\) and \( N \) and are polynomials in \( \omega \) of degree \( p \) with coefficients that are polynomials in \( t \) of degree at most \( 2p \), while the functions \( r_{n,N}(t) \) are bounded in absolute value for \( 0 \leq t \leq a \) by a polynomial of degree \( 2N \) whose coefficients are independent of \(n\). Finally, we get from Lemmas~\ref{lem:7} and~\ref{lem:8} that
\[
\sum_{j=0}^1 \beta_j(t)\alpha_{1-j}(t) = t+ \mathcal O\big(t^2\big) \quad \text{and} \quad \sum_{j=0}^k \beta_j(t)\alpha_{k-j}(t) = \mathcal O\big(t^2\big)
\]
for all \( k\geq 2 \). Therefore, it holds that \( r_p(t) = \mathcal O\big(t^2\big) \) as \( t\to 0 \) for all \( p\geq 1\).

It follows from Lemma~\ref{lem:6} that \(E_n(t) = t^{-2} \chi(t) [1-R_n(t)] \). Hence, plugging the expansion of \( R_n(t) \) into this formula gives us
\[
E_n(t)= \chi(t) \left[ \sum_{p=1}^{N-1} e_{p}(t)(n+1)^{-p} + e_{n,N}(t)(n+1)^{-N}\right],
\]
where \( e_p(t) := - t^{-2} r_p(t) \) for any \( p \) and  \( e_{n,N}(t) := - t^{-2} r_{n,N}(t) \) for any \( n,N \). It follows from the properties of \( r_p(t) \) that each \( e_p(t) \) is a continuous function and is bounded in absolute value on \( 0\leq t < \infty \) by a polynomial of degree \( 2p-2 \). Also, since \( \chi(t) \) is a continuous function as well and \( \lim_{t\to0^+}E_n(t) \) exists and is finite according to Lemma~\ref{lem:6}, so must \( \lim_{t\to 0^+}e_{n,N}(t) \) for all \( {n,N} \). Then it follows from properties of \( r_{n,N}(t)\) that \( e_{n,N}(t) \) is bounded in absolute value when \( 0 \leq t \leq a \) by a polynomial of degree \( 2N-2 \) whose coefficients are independent of \(n\). 

From what precedes, we get that
\[
|E_n(t)| \leq \frac{\chi(t)|e_{n,1}(t)|}{n+1} \leq \frac{C_{\mu,T}}{n+1}, \quad 0\leq t\leq a.
\]
Hence, for any \( -1<E<0 \) there exists an integer \( n_E \) such that \( E\leq E_n(t) \) for all \( 0\leq t \leq a \) and \( n\geq n_E \). Thus, by applying \eqref{taylor} one more time, we get that
\[
(1+E_n(t))^{1/2}-1 = \sum_{k=1}^{N-1} \binom{1/2}{k}E_n^k(t) + \binom{1/2}{N} \frac{E_n^N(t)}{(1+\theta_5 E_n(t))^{N-1/2}}
\]
for some \( \theta_5 \in (0,1) \) that depends on \( N \) and \( E_n(t) \). Therefore,  the claim of the lemma follows with
\[
u_{p}(t) := \sum_{k=1}^p \binom{1/2}{k} \chi^{k-1}(t) \sum_{j_1+\cdots+j_k=p} e_{j_1}(t) \cdots e_{j_k}(t),
\]
which is bounded in absolute value on \( 0\leq t <\infty \) by a polynomial of degree $ 2p-2 $ whose coefficients are independent of $n$ and $N$, and
\[
u_{n,N}(t) := \sum_{k=1}^{N-1}  \binom{1/2}{k} \chi^{k-1}(t) \sum_{j_1+\cdots+j_k\geq N} \frac{e_{n,j_1,N}(t) \cdots e_{n,j_k,N}(t)}{(n+1)^{j_1+\cdots+j_k-N}}+ \binom{1/2}{N} \frac{(n+1)^NE_n^N(t)}{(1+\theta_5 E_n(t))^{N-1/2}}
\]
where \( e_{n,j,N}(t) := e_j(t) \) when \( j<N \) and \( e_{n,N,N}(t) := e_{n,N}(t) \), which is bounded in absolute value on \( 0 \leq t \leq a \) by a polynomial of degree $ 2N-2 $ whose coefficients are independent of $n$ due to the same reasoning as in two previous lemmas.
\qed\endproof

\begin{lemma}
\label{lem:10}
Given \( N\geq 1 \), it holds that
\[
\frac{(1+E_n(t))^{1/2}-1}{2(n+1)-t} = \chi(t)\sum_{p=2}^{N-1}v_p(t)(n+1)^{-p} + \chi(t)v_{n,N}(t)(n+1)^{-N},
\]
where \( v_p(t) \) is bounded in absolute value on \( 0\leq t< \infty \) by a polynomial of degree \( 2p-4 \) whose coefficients are independent of $n$ and $N$ and the functions \( v_{n,N}(t) \) is bounded in absolute value when \( 0 \leq t \leq a \) by a polynomial of degree $ 2N-4$ whose coefficients are independent of $n$.
\end{lemma}
\proof
Since \( 0 \leq t\leq a=\sqrt{n+1} \), we get from \eqref{taylor} that
\[
\frac1{2(n+1)-t} = \sum_{p=1}^{N-1}z_p(t)(n+1)^{-p} + z_{n,N}(t)(n+1)^{-N},
\]
where
\[
z_p(t) := 2^{-p}t^{p-1} \quad \text{and} \quad z_{n,N}(t) :=  \frac{2^{-N}t^{N-1}}{(1-\theta_6 t/2(n+1))^{N+1}}
\]
for some \( \theta_6\in(0,1) \) that depends on \( N \) and \( t \). Therefore, the claim of the lemma follows from Lemma~\ref{lem:9} with
\[
v_p(t) := \sum_{j=1}^{p-1}z_j(t)u_{p-j}(t) \quad \text{and} \quad v_{n,N}(t) := \sum_{k=N}^{2N}\sum_{j_1+j_2=k}\frac{z_{n,j_1,N}(t) v_{n,j_2,N}(t)}{(n+1)^{k-N}}
\]
where \( j_1,j_2\in\{1,\ldots,N\} \), \( z_{n,j,N}(t) := z_j(t) \), \( u_{n,j,N}(t) := u_j(t) \) for \( j<N \), and \( z_{n,n,N}(t) := z_{n,N}(t) \), \( u_{n,N,N}(t) := u_{n,N}(t) \).
\qed\endproof

With the notation introduced in Lemmas~\ref{lem:5},~\ref{lem:9}, and~\ref{lem:10}, the following lemma holds.

\begin{lemma}
\label{lem:11}
Given \( N\geq 1 \), it holds that
\[
G_n(t) = I_1^\mu (n+1)^{-1} + \sum_{p=2}^{N-1}\big( I_p^\mu + J_p^\mu \big)(n+1)^{-p} + \mathcal O_N \left( (n+1)^{-N} \right)
\]
for all \( n \) large, where 
\[
I_p^\mu := \frac1\pi \int_{0}^{\infty} t^{-1}f(t)\chi(t)u_p(t)\dd t \quad \text{and} \quad
J_p^\mu := \frac1\pi \int_{0}^{\infty} f(t)\chi(t)v_p(t) \dd t
\]
(observe that \( t^{-1}f(t) \) is a continuous and bounded function on \( 0\leq t<\infty \), \( \chi(t) \) decreases exponentially at infinity, and the functions \( u_p(t),v_p(t) \) are bounded by polynomials).
\end{lemma}
\proof
By the very definition of \( G_n(t) \) in Lemma~\ref{lem:5} we have that \( G_n(t) = I_n(t) + J_n(t) \), where
\[
I_n(t) := \frac1\pi\int_0^a t^{-1}f(t)\left((1+E_n(t))^{1/2}-1\right)\dd t
\]
and
\[
J_n(t) := \frac1\pi\int_0^a f(t)\frac{(1+E_n(t))^{1/2}-1}{2(n+1)-t}\dd t.
\]
Using Lemma~\ref{lem:9}, we can rewrite the first integral above as
\[
I_n(t) = \sum_{p=1}^{N-1} I_p^\mu(n+1)^{-p} -S_n(t) + T_n(t), 
\]
where
\[
S_n(t) := \frac{1}{\pi} \sum_{p=1}^{N-1} (n+1)^{-p} \int_{a}^{\infty} t^{-1}f(t)\chi(t)u_p(t)\dd t
\]
and
\[
T_n(t) :=  \frac{1}{\pi} (n+1)^{-N} \int_{0}^{a} t^{-1}f(t)\chi(t) u_{n,N}(t)\dd t.
\]
Since \( u_p(t) = \mathcal{O}\big(t^{2p-2}\big) \), \( f(t)=\mathcal O(1) \), and \( \chi(t) = \mathcal O\big(t^4e^{-2t}\big) \) as \( t\to\infty \), it holds that
\begin{multline*}
S_n(t) = \sum_{p=1}^{N-1} (n+1)^{-p}  \int_{a}^{\infty} \mathcal{O}\big(t^{2p+1} e^{-2t}\big) \dd t = \sum_{p=1}^{N-1} (n+1)^{-p} \mathcal{O}\big(a^{2p+1}e^{-2a}\big) = \\ = \mathcal{O}_N\left(ae^{-2a}\right) = o_N\left((n+1)^{-N}\right).
\end{multline*}
Moreover, since \( u_{n,N}(t) \) is bounded by a polynomial of degree \( 2N-2 \) for \( 0\leq t\leq a \), we have that \( T_n(t) = \mathcal{O}_N\big( (n+1)^{-N} \big)\).

Similarly, we get from Lemma~\ref{lem:10} that
\[
J_n(t) = \sum_{p=2}^{N-1}J_{p}^\mu(n+1)^{-p} - U_n(t)+ V_n(t),
\]
where
\[
U_n(t) := \frac1\pi \sum_{p=2}^{N-1}(n+1)^{-p} \int_{a}^{\infty} f(t)\chi(t)v_{p}(t) \dd t
\]
and
\[
V_n(t) := \frac1\pi(n+1)^{-N} \int_{0}^{a} f(t)\chi(t) v_{n,N}(t) \dd t.
\]
An argument as above argument shows that \( U_n(t) = \mathcal{O}_N\big(e^{-2a}\big)=o_N\big( (n+1)^{-N} \big) \) and \( V_n(t) = \mathcal{O}_N\big( (n+1)^{-N} \big) \) for large $ n $, which finishes the proof of the lemma.
\qed\endproof

\begin{lemma}
The claim of Theorem~\ref{thm:1} holds.
\end{lemma}
\proof
It follows from Lemmas~\ref{lem:5} and~\ref{lem:11} that given an integer \( N\geq 1 \), it holds that
\[
\widehat{\mathbb E}_n(\mu) = \frac1\pi\log(n+1) + \frac12A_0  + \sum_{p=1}^{N-1}\big(I_p^\mu + J_p^\mu - H_p/2\big)(n+1)^{-p}  + \mathcal O_N\left((n+1)^{-N}\right),
\]
where we set \( J_1^\mu:=0 \). The claim of Theorem~\ref{thm:1} now follows from \eqref{reduction} by taking \( A_p^\mu := I_p^\mu + I_p^\sigma + J_p^\mu + J_p^\sigma - H_p \).
\qed\endproof

\section*{Acknowledgments}

The work of the first author is done towards completion of her Ph.D. degree at Indiana University-Purdue University Indianapolis under the direction of the second author. The research of the second author was supported in part by a grant from the Simons Foundation, CGM-354538.

\small

\bibliographystyle{elsarticle-num}

\end{document}